\documentclass{amsart}
\usepackage{amsfonts}

\newtheorem{theorem}{Theorem}
\theoremstyle{plain}

\newtheorem{corollary}{Corollary}

\newtheorem{definition}{Definition}

\newtheorem{lemma}{Lemma}

\newtheorem{proposition}{Proposition}
\newtheorem{remark}{Remark}

\numberwithin{equation}{section}
\input{tcilatex}

\begin{document}
\title[Additive Reverses of the Generalised Triangle Inequality]{Additive
Reverses of the Generalised Triangle Inequality in Normed Spaces}
\author{S.S. Dragomir}
\address{School of Computer Science and Mathematics\\
Victoria University of Technology\\
PO Box 14428, MCMC 8001\\
Victoria, Australia.}
\email{sever.dragomir@vu.edu.au}
\urladdr{http://rgmia.vu.edu.au/SSDragomirWeb.html}
\date{June 21, 2004.}

\begin{abstract}
Some additive reverses of the generalised triangle inequality in normed
linear spaces are given. Applications for complex numbers are provided as
well.
\end{abstract}

\maketitle

\section{Introduction}

In \cite{DM}, Diaz and Metcalf established the following reverse of the
generalised triangle inequality in real or complex normed linear spaces.

If $F:X\rightarrow \mathbb{K}$, $\mathbb{K}=\mathbb{R},\mathbb{C}$ is a
linear functional of a unit norm defined on the normed linear space $X$
endowed with the norm $\left\Vert \cdot \right\Vert $ and the vectors $%
x_{1},\dots ,x_{n}$ satisfy the condition%
\begin{equation}
0\leq r\left\Vert x_{i}\right\Vert \leq \func{Re}F\left( x_{i}\right) ,\ \ \
\ \ i\in \left\{ 1,\dots ,n\right\} ;  \label{1.1}
\end{equation}%
then%
\begin{equation}
r\sum_{i=1}^{n}\left\Vert x_{i}\right\Vert \leq \left\Vert
\sum_{i=1}^{n}x_{i}\right\Vert ,  \label{1.2}
\end{equation}%
where equality holds if and only if both%
\begin{equation}
F\left( \sum_{i=1}^{n}x_{i}\right) =r\sum_{i=1}^{n}\left\Vert
x_{i}\right\Vert  \label{1.3}
\end{equation}%
and%
\begin{equation}
F\left( \sum_{i=1}^{n}x_{i}\right) =\left\Vert
\sum_{i=1}^{n}x_{i}\right\Vert .  \label{1.4}
\end{equation}

If $X=H$, $\left( H;\left\langle \cdot ,\cdot \right\rangle \right) $ is an
inner product space and $F\left( x\right) =\left\langle x,e\right\rangle ,$ $%
\left\Vert e\right\Vert =1,$ then the condition (\ref{1.1}) may be replaced
with the simpler assumption%
\begin{equation}
0\leq r\left\Vert x_{i}\right\Vert \leq \func{Re}\left\langle
x_{i},e\right\rangle ,\qquad i=1,\dots ,n,  \label{1.5}
\end{equation}%
which implies the reverse of the generalised triangle inequality (\ref{1.2}%
). In this case the equality holds in (\ref{1.2}) if and only if \cite{DM}%
\begin{equation}
\sum_{i=1}^{n}x_{i}=r\left( \sum_{i=1}^{n}\left\Vert x_{i}\right\Vert
\right) e.  \label{1.6}
\end{equation}%
Let $F_{1},\dots ,F_{m}$ be linear functionals on $X,$ each of unit norm.
Let \cite{DM}%
\begin{equation*}
c=\sup_{x\neq 0}\left[ \frac{\sum_{k=1}^{m}\left\vert F_{k}\left( x\right)
\right\vert ^{2}}{\left\Vert x\right\Vert ^{2}}\right] ;
\end{equation*}%
it then follows that $1\leq c\leq m.$ Suppose the vectors $x_{1},\dots ,x_{k}
$ whenever $x_{i}\neq 0,$ satisfy%
\begin{equation}
0\leq r_{k}\left\Vert x_{i}\right\Vert \leq \func{Re}F_{k}\left(
x_{i}\right) ,\qquad i=1,\dots ,n,\ k=1,\dots ,m.  \label{1.7}
\end{equation}%
Then \cite{DM}%
\begin{equation}
\left( \frac{\sum_{k=1}^{m}r_{k}^{2}}{c}\right) \sum_{i=1}^{n}\left\Vert
x_{i}\right\Vert \leq \left\Vert \sum_{i=1}^{n}x_{i}\right\Vert ,
\label{1.8}
\end{equation}%
where equality holds if and only if both%
\begin{equation}
F_{k}\left( \sum_{i=1}^{n}x_{i}\right) =r_{k}\sum_{i=1}^{n}\left\Vert
x_{i}\right\Vert ,\qquad k=1,\dots ,m  \label{1.9}
\end{equation}%
and 
\begin{equation}
\sum_{k=1}^{m}\left[ F_{k}\left( \sum_{i=1}^{n}x_{i}\right) \right]
^{2}=c\left\Vert \sum_{i=1}^{n}x_{i}\right\Vert ^{2}.  \label{1.10}
\end{equation}

If $X=H,$ an inner product space, then, for $F_{k}\left( x\right)
=\left\langle x,e_{k}\right\rangle ,$ where $\left\{ e_{k}\right\} _{k=%
\overline{1,n}}$ is an orthonormal family in $H,$ i.e., $\left\langle
e_{i},e_{j}\right\rangle =\delta _{ij},$ $i,j\in \left\{ 1,\dots ,k\right\} ,
$ $\delta _{ij}$ is Kronecker delta, the condition (\ref{1.7}) may be
replaced by%
\begin{equation}
0\leq r_{k}\left\Vert x_{i}\right\Vert \leq \func{Re}\left\langle
x_{i},e\right\rangle ,\qquad i=1,\dots ,n,\ k=1,\dots ,m.  \label{1.11}
\end{equation}%
implying the following reverse of the generalised triangle inequality%
\begin{equation}
\left( \sum_{k=1}^{m}r_{k}^{2}\right) ^{\frac{1}{2}}\sum_{i=1}^{n}\left\Vert
x_{i}\right\Vert \leq \left\Vert \sum_{i=1}^{n}x_{i}\right\Vert ,
\label{1.12}
\end{equation}%
where the equality holds if and only if%
\begin{equation}
\sum_{i=1}^{n}x_{i}=\left( \sum_{i=1}^{n}\left\Vert x_{i}\right\Vert \right)
\sum_{k=1}^{m}r_{k}e_{k}.  \label{1.13}
\end{equation}

The main aim of this paper is to provide some new reverse results of the
generalised triangle inequality in its additive form, namely, upper bounds
for the nonnegative quantity%
\begin{equation*}
\sum_{i=1}^{n}\left\Vert x_{i}\right\Vert -\left\Vert
\sum_{i=1}^{n}x_{i}\right\Vert
\end{equation*}%
under various assumptions for the vectors $x_{i},$ $i\in \left\{ 1,\dots
,n\right\} $ in a real or complex normed space $\left( X,\left\Vert \cdot
\right\Vert \right) .$ Applications for complex numbers are provided as well.

\section{Semi-Inner Products and Diaz-Metcalf Inequality}

In 1961, G. Lumer \cite{L} introduced the following concept.

\begin{definition}
\label{d3.1}Let $X$ be a linear space over the real or complex number field $%
\mathbb{K}$. The mapping $\left[ \cdot ,\cdot \right] :X\times X\rightarrow 
\mathbb{K}$ is called a \textit{semi-inner product }on $X,$ if the following
properties are satisfied (see also \cite[p. 17]{SSD1}):

\begin{enumerate}
\item[$\left( i\right) $] $\left[ x+y,z\right] =\left[ x,z\right] +\left[ y,z%
\right] $ for all $x,y,z\in X;$

\item[$\left( ii\right) $] $\left[ \lambda x,y\right] =\lambda \left[ x,y%
\right] $ for all $x,y\in X$ and $\lambda \in \mathbb{K}$;

\item[$\left( iii\right) $] $\left[ x,x\right] \geq 0$ for all $x\in X$ and $%
\left[ x,x\right] =0$ implies $x=0$;

\item[$\left( iv\right) $] $\left\vert \left[ x,y\right] \right\vert
^{2}\leq \left[ x,x\right] \left[ y,y\right] $ for all $x,y\in X;$

\item[$\left( v\right) $] $\left[ x,\lambda y\right] =\bar{\lambda}\left[ x,y%
\right] $ for all $x,y\in X$ and $\lambda \in \mathbb{K}$.
\end{enumerate}
\end{definition}

It is well known that the mapping $X\ni x\longmapsto \left[ x,x\right] ^{%
\frac{1}{2}}\in \mathbb{R}$ is a norm on $X$ and for any $y\in X,$ the
functional $X\ni x\overset{\varphi _{y}}{\longmapsto }\left[ x,x\right] ^{%
\frac{1}{2}}\in \mathbb{K}$ is a continuous linear functional on $X$ endowed
with the norm $\left\Vert \cdot \right\Vert $ generated by $\left[ \cdot
,\cdot \right] .$ Moreover, one has $\left\Vert \varphi _{y}\right\Vert
=\left\Vert y\right\Vert $ (see for instance \cite[p. 17]{SSD1}).

Let $\left( X,\left\Vert \cdot \right\Vert \right) $ be a real or complex
normed space. If $J:X\rightarrow _{2}X^{\ast }$ is the \textit{normalised
duality mapping }defined on $X,$ i.e., we recall that (see for instance \cite%
[p. 1]{SSD1})%
\begin{equation*}
J\left( X\right) =\left\{ \varphi \in X^{\ast }|\varphi \left( x\right)
=\left\Vert \varphi \right\Vert \left\Vert x\right\Vert ,\ \left\Vert
\varphi \right\Vert =\left\Vert x\right\Vert \right\} ,\ \ \ x\in X,
\end{equation*}%
then we may state the following representation result (see for instance \cite%
[p. 18]{SSD1}):

Each semi-inner product $\left[ \cdot ,\cdot \right] :X\times X\rightarrow 
\mathbb{K}$ that generates the norm $\left\Vert \cdot \right\Vert $ of the
normed linear space $\left( X,\left\Vert \cdot \right\Vert \right) $ over
the real or complex number field $\mathbb{K}$, is of the form%
\begin{equation*}
\left[ x,y\right] =\left\langle \tilde{J}\left( y\right) ,x\right\rangle 
\text{ \ for any \ }x,y\in X,
\end{equation*}%
where $\tilde{J}$ is a selection of the normalised duality mapping and $%
\left\langle \varphi ,x\right\rangle :=\varphi \left( x\right) $ for $%
\varphi \in X^{\ast }$ and $x\in X.$

Utilising the concept of semi-inner products, we can state the following
particular case of the Diaz-Metcalf inequality.

\begin{corollary}
\label{c3.1}Let $\left( X,\left\Vert \cdot \right\Vert \right) $ be a normed
linear space, $\left[ \cdot ,\cdot \right] :X\times X\rightarrow \mathbb{K}$
a semi-inner product generating the norm $\left\Vert \cdot \right\Vert $ and 
$e\in X,$ $\left\Vert e\right\Vert =1.$ If $x_{i}\in X,$ $i\in \left\{
1,\dots ,n\right\} $ and $r\geq 0$ such that%
\begin{equation}
r\left\Vert x_{i}\right\Vert \leq \func{Re}\left[ x_{i},e\right] \text{ \
for each \ }i\in \left\{ 1,\dots ,n\right\} ,  \label{3.1}
\end{equation}%
then we have the inequality%
\begin{equation}
r\sum_{i=1}^{n}\left\Vert x_{i}\right\Vert =\left\Vert
\sum_{i=1}^{n}x_{i}\right\Vert .  \label{3.2}
\end{equation}%
The case of equality holds in (\ref{3.2}) if and only if both%
\begin{equation}
\left[ \sum_{i=1}^{n}x_{i},e\right] =r\sum_{i=1}^{n}\left\Vert
x_{i}\right\Vert  \label{3.3}
\end{equation}%
and 
\begin{equation}
\left[ \sum_{i=1}^{n}x_{i},e\right] =\left\Vert
\sum_{i=1}^{n}x_{i}\right\Vert .  \label{3.4}
\end{equation}
\end{corollary}

The proof is obvious from the Diaz-Metcalf theorem \cite[Theorem 3]{DM}
applied for the continuous linear functional $F_{e}\left( x\right) =\left[
x,e\right] ,$ $x\in X.$

Before we provide a simpler necessary and sufficient condition of equality
in (\ref{3.2}), we need to recall the concept of strictly convex normed
spaces and a classical characterisation of these spaces.

\begin{definition}
\label{d3.2}A normed linear space $\left( X,\left\Vert \cdot \right\Vert
\right) $ is said to be strictly convex if for every $x,y$ from $X$ with $%
x\neq y$ and $\left\Vert x\right\Vert =\left\Vert y\right\Vert =1,$ we have $%
\left\Vert \lambda x+\left( 1-\lambda \right) y\right\Vert <1$ for all $%
\lambda \in \left( 0,1\right) .$
\end{definition}

The following characterisation of strictly convex spaces is useful in what
follows (see \cite{B}, \cite{GS}, or \cite[p. 21]{SSD1}).

\begin{theorem}
\label{t3.1}Let $\left( X,\left\Vert \cdot \right\Vert \right) $ be a normed
linear space over $\mathbb{K}$ and $\left[ \cdot ,\cdot \right] $ a
semi-inner product generating its norm. The following statements are
equivalent:

\begin{enumerate}
\item[$\left( i\right) $] $\left( X,\left\Vert \cdot \right\Vert \right) $
is strictly convex;

\item[$\left( ii\right) $] For every $x,y\in X,$ $x,y\neq 0$ with $\left[ x,y%
\right] =\left\Vert x\right\Vert \left\Vert y\right\Vert ,$ there exists a $%
\lambda >0$ such that $x=\lambda y.$
\end{enumerate}
\end{theorem}

The following result may be stated.

\begin{corollary}
\label{c3.2}Let $\left( X,\left\Vert \cdot \right\Vert \right) $ be a
strictly convex normed linear space, $\left[ \cdot ,\cdot \right] $ a
semi-inner product generating the norm and $e,$ $x_{i}$ $\left( i\in \left\{
1,\dots ,n\right\} \right) $ as in Corollary \ref{c3.1}. Then the case of
equality holds in (\ref{3.2}) if and only if%
\begin{equation}
\sum_{i=1}^{n}x_{i}=r\left( \sum_{i=1}^{n}\left\Vert x_{i}\right\Vert
\right) e.  \label{3.5}
\end{equation}
\end{corollary}

\begin{proof}
If (\ref{3.5}) holds true, then, obviously%
\begin{equation*}
\left\Vert \sum_{i=1}^{n}x_{i}\right\Vert =r\left( \sum_{i=1}^{n}\left\Vert
x_{i}\right\Vert \right) \left\Vert e\right\Vert =r\sum_{i=1}^{n}\left\Vert
x_{i}\right\Vert ,
\end{equation*}%
which is the equality case in (\ref{3.2}).

Conversely, if the equality holds in (\ref{3.2}), then by Corollary \ref%
{c3.1}, we have that (\ref{3.3}) and (\ref{3.4}) hold true. Utilising
Theorem \ref{t3.1}, we conclude that there exists a $\mu >0$ such that%
\begin{equation}
\sum_{i=1}^{n}x_{i}=\mu e.  \label{3.6}
\end{equation}%
Inserting this in (\ref{3.3}) we get%
\begin{equation*}
\mu \left\Vert e\right\Vert ^{2}=r\sum_{i=1}^{n}\left\Vert x_{i}\right\Vert
\end{equation*}%
giving%
\begin{equation}
\mu =r\sum_{i=1}^{n}\left\Vert x_{i}\right\Vert .  \label{3.7}
\end{equation}%
Finally, by (\ref{3.6}) and (\ref{3.7}) we deduce (\ref{3.5}) and the
corollary is proved.
\end{proof}

\section{An Additive Reverse for the Triangle Inequality}

In the following we provide an alternative of the Diaz-Metcalf reverse of
the generalised triangle inequality.

\begin{theorem}
\label{t5.1}Let $\left( X,\left\Vert \cdot \right\Vert \right) $ be a normed
linear space over the real or complex number field $\mathbb{K}$ and $%
F:X\rightarrow \mathbb{K}$ a linear functional with the property that $%
\left\vert F\left( x\right) \right\vert \leq \left\Vert x\right\Vert $ for
any $x\in X$ (i.e., $\left\Vert F\right\Vert =1,$ we say that $F$ is of unit
norm). If $x_{i}\in X,$ $k_{i}\geq 0,$ $i\in \left\{ 1,\dots ,n\right\} $
are such that%
\begin{equation}
\left( 0\leq \right) \left\Vert x_{i}\right\Vert -\func{Re}F\left(
x_{i}\right) \leq k_{i}\text{ \ for each \ }i\in \left\{ 1,\dots ,n\right\} ,
\label{5.1}
\end{equation}%
then we have the inequality%
\begin{equation}
\left( 0\leq \right) \sum_{i=1}^{n}\left\Vert x_{i}\right\Vert -\left\Vert
\sum_{i=1}^{n}x_{i}\right\Vert \leq \sum_{i=1}^{n}k_{i}.  \label{5.2}
\end{equation}%
The equality holds in (\ref{5.2}) if and only if both%
\begin{equation}
F\left( \sum_{i=1}^{n}x_{i}\right) =\left\Vert
\sum_{i=1}^{n}x_{i}\right\Vert \text{ \ and \ }F\left(
\sum_{i=1}^{n}x_{i}\right) =\sum_{i=1}^{n}\left\Vert x_{i}\right\Vert
-\sum_{i=1}^{n}k_{i}.  \label{5.3}
\end{equation}
\end{theorem}

\begin{proof}
If we sum in (\ref{5.1}) over $i$ from $1$ to $n,$ then we get%
\begin{equation}
\sum_{i=1}^{n}\left\Vert x_{i}\right\Vert \leq \func{Re}\left[ F\left(
\sum_{i=1}^{n}x_{i}\right) \right] +\sum_{i=1}^{n}k_{i}.  \label{5.4}
\end{equation}%
Taking into account that $\left\vert F\left( x\right) \right\vert \leq
\left\Vert x\right\Vert $ for each $x\in X,$ then we may state that%
\begin{equation}
\func{Re}\left[ F\left( \sum_{i=1}^{n}x_{i}\right) \right] \leq \left\vert 
\func{Re}F\left( \sum_{i=1}^{n}x_{i}\right) \right\vert \leq \left\vert
F\left( \sum_{i=1}^{n}x_{i}\right) \right\vert \leq \left\Vert
\sum_{i=1}^{n}x_{i}\right\Vert .  \label{5.5}
\end{equation}%
Now, making use of (\ref{5.4}) and (\ref{5.5}), we deduce (\ref{5.2}).

Obviously, if (\ref{5.3}) is valid, then the case of equality in (\ref{5.2})
holds true.

Conversely, if the equality holds in (\ref{5.2}), then it must hold in all
the inequalities used to prove (\ref{5.2}), therefore we have%
\begin{equation*}
\sum_{i=1}^{n}\left\Vert x_{i}\right\Vert =\func{Re}\left[ F\left(
\sum_{i=1}^{n}x_{i}\right) \right] +\sum_{i=1}^{n}k_{i}
\end{equation*}%
and%
\begin{equation*}
\func{Re}\left[ F\left( \sum_{i=1}^{n}x_{i}\right) \right] =\left\vert
F\left( \sum_{i=1}^{n}x_{i}\right) \right\vert =\left\Vert
\sum_{i=1}^{n}x_{i}\right\Vert ,
\end{equation*}%
which imply (\ref{5.3}).
\end{proof}

The following corollary may be stated.

\begin{corollary}
\label{c5.1}Let $\left( X,\left\Vert \cdot \right\Vert \right) $ be a normed
linear space, $\left[ \cdot ,\cdot \right] :X\times X\rightarrow \mathbb{K}$
a semi-inner product generating the norm $\left\Vert \cdot \right\Vert $ and 
$e\in X,$ $\left\Vert e\right\Vert =1.$ If $x_{i}\in X,$ $k_{i}\geq 0,$ $\
i\in \left\{ 1,\dots ,n\right\} $ are such that%
\begin{equation}
\left( 0\leq \right) \left\Vert x_{i}\right\Vert -\func{Re}\left[ x_{i},e%
\right] \leq k_{i}\text{ \ for each \ }i\in \left\{ 1,\dots ,n\right\} ,
\label{5.6}
\end{equation}%
then we have the inequality%
\begin{equation}
\left( 0\leq \right) \sum_{i=1}^{n}\left\Vert x_{i}\right\Vert -\left\Vert
\sum_{i=1}^{n}x_{i}\right\Vert \leq \sum_{i=1}^{n}k_{i}.  \label{5.7}
\end{equation}%
The equality holds in (\ref{5.7}) if and only if both%
\begin{equation}
\left[ \sum_{i=1}^{n}x_{i},e\right] =\left\Vert
\sum_{i=1}^{n}x_{i}\right\Vert \text{ \ \ and \ \ }\left[
\sum_{i=1}^{n}x_{i},e\right] =\sum_{i=1}^{n}\left\Vert x_{i}\right\Vert
-\sum_{i=1}^{n}k_{i}.  \label{5.8}
\end{equation}%
Moreover, if $\left( X,\left\Vert \cdot \right\Vert \right) $ is strictly
convex, then the case of equality holds in (\ref{5.7}) if and only if%
\begin{equation}
\sum_{i=1}^{n}\left\Vert x_{i}\right\Vert \geq \sum_{i=1}^{n}k_{i}
\label{5.9}
\end{equation}%
and 
\begin{equation}
\sum_{i=1}^{n}x_{i}=\left( \sum_{i=1}^{n}\left\Vert x_{i}\right\Vert
-\sum_{i=1}^{n}k_{i}\right) \cdot e.  \label{5.10}
\end{equation}
\end{corollary}

\begin{proof}
The first part of the corollary is obvious by Theorem \ref{t5.1} applied for
the continuous linear functional of unit norm $F_{e},$ $F_{e}\left( x\right)
=\left[ x,e\right] ,$ $x\in X.$ The second part may be shown on utilising a
similar argument to the one from the proof of Corollary \ref{c3.2}. We omit
the details.
\end{proof}

\begin{remark}
\label{r5.1}If $X=H,$ $\left( H;\left\langle \cdot ,\cdot \right\rangle
\right) $ is an inner product space, then from Corollary \ref{c5.1} we
deduce the additive reverse inequality obtained in Theorem 7 of \cite{SSD2}.
For further similar results in inner product spaces, see \cite{SSD2} and 
\cite{SSD3}.
\end{remark}

\section{Reverse Inequalities for $m$ Functionals}

The following result generalising Theorem \ref{t5.1} may be stated.

\begin{theorem}
\label{t2.1}Let $\left( X,\left\Vert \cdot \right\Vert \right) $ be a normed
linear space over the real or complex number field $\mathbb{K}$. If $F_{k}$, 
$k\in \left\{ 1,\dots ,m\right\} $ are bounded linear functionals defined on 
$X\ $and $x_{i}\in X,$ $M_{ik}\geq 0$ for $i\in \left\{ 1,\dots ,n\right\} $%
, $k\in \left\{ 1,\dots ,m\right\} $ such that%
\begin{equation}
\left\Vert x_{i}\right\Vert -\func{Re}F_{k}\left( x_{i}\right) \leq M_{ik}%
\text{ \ for each \ }i\in \left\{ 1,\dots ,n\right\} ,\ k\in \left\{ 1,\dots
,m\right\} ,  \label{2.1}
\end{equation}%
then we have the inequality%
\begin{equation}
\sum_{i=1}^{n}\left\Vert x_{i}\right\Vert \leq \left\Vert \frac{1}{m}%
\sum_{k=1}^{m}F_{k}\right\Vert \left\Vert \sum_{i=1}^{n}x_{i}\right\Vert +%
\frac{1}{m}\sum_{k=1}^{m}\sum_{i=1}^{n}M_{ik}.  \label{2.2}
\end{equation}%
The case of equality holds in (\ref{2.2}) if both%
\begin{equation}
\frac{1}{m}\sum_{k=1}^{m}F_{k}\left( \sum_{i=1}^{n}x_{i}\right) =\left\Vert 
\frac{1}{m}\sum_{k=1}^{m}F_{k}\right\Vert \left\Vert
\sum_{i=1}^{n}x_{i}\right\Vert  \label{2.3}
\end{equation}%
and%
\begin{equation}
\frac{1}{m}\sum_{k=1}^{m}F_{k}\left( \sum_{i=1}^{n}x_{i}\right)
=\sum_{i=1}^{n}\left\Vert x_{i}\right\Vert -\frac{1}{m}\sum_{k=1}^{m}%
\sum_{j=1}^{n}M_{ik}.  \label{2.4}
\end{equation}
\end{theorem}

\begin{proof}
If we sum (\ref{2.1}) over $i$ from $1$ to $n,$ then we deduce%
\begin{equation*}
\sum_{i=1}^{n}\left\Vert x_{i}\right\Vert -\func{Re}F_{k}\left(
\sum_{i=1}^{n}x_{i}\right) \leq \sum_{i=1}^{n}M_{ik}
\end{equation*}%
for each $k\in \left\{ 1,\dots ,m\right\} .$

Summing these inequalities over $k$ from $1$ to $m,$ we deduce 
\begin{equation}
\sum_{i=1}^{n}\left\Vert x_{i}\right\Vert \leq \frac{1}{m}\sum_{k=1}^{m}%
\func{Re}F_{k}\left( \sum_{i=1}^{n}x_{i}\right) +\frac{1}{m}%
\sum_{k=1}^{m}\sum_{i=1}^{n}M_{ik}.  \label{2.5}
\end{equation}%
Utilising the continuity property of the functionals $F_{k}$ and the
properties of the modulus, we have%
\begin{align}
\sum_{k=1}^{m}\func{Re}F_{k}\left( \sum_{i=1}^{n}x_{i}\right) & \leq
\left\vert \sum_{k=1}^{m}\func{Re}F_{k}\left( \sum_{i=1}^{n}x_{i}\right)
\right\vert  \label{2.6} \\
& \leq \left\vert \sum_{k=1}^{m}F_{k}\left( \sum_{i=1}^{n}x_{i}\right)
\right\vert \leq \left\Vert \sum_{k=1}^{m}F_{k}\right\Vert \left\Vert
\sum_{i=1}^{n}x_{i}\right\Vert .  \notag
\end{align}%
Now, by (\ref{2.5}) and (\ref{2.6}), we deduce (\ref{2.2}).

Obviously, if (\ref{2.3}) and (\ref{2.4}) hold true, then the case of
equality is valid in (\ref{2.2}).

Conversely, if the case of equality holds in (\ref{2.2}), then it must hold
in all the inequalities used to prove (\ref{2.2}). Therefore we have%
\begin{equation*}
\sum_{i=1}^{n}\left\Vert x_{i}\right\Vert =\frac{1}{m}\sum_{k=1}^{m}\func{Re}%
F_{k}\left( \sum_{i=1}^{n}x_{i}\right) +\frac{1}{m}\sum_{k=1}^{m}%
\sum_{i=1}^{n}M_{ik},
\end{equation*}%
\begin{equation*}
\sum_{k=1}^{m}\func{Re}F_{k}\left( \sum_{i=1}^{n}x_{i}\right) =\left\Vert
\sum_{k=1}^{m}F_{k}\right\Vert \left\Vert \sum_{i=1}^{n}x_{i}\right\Vert
\end{equation*}%
and%
\begin{equation*}
\sum_{k=1}^{m}\func{Im}F_{k}\left( \sum_{i=1}^{n}x_{i}\right) =0.
\end{equation*}%
These imply that (\ref{2.3}) and (\ref{2.4}) hold true, and the theorem is
completely proved.
\end{proof}

\begin{remark}
\label{r2.1}If $F_{k},$ $k\in \left\{ 1,\dots ,m\right\} $ are of unit norm,
then, from (\ref{2.2}), we deduce the inequality%
\begin{equation}
\sum_{i=1}^{n}\left\Vert x_{i}\right\Vert \leq \left\Vert
\sum_{i=1}^{n}x_{i}\right\Vert +\frac{1}{m}\sum_{k=1}^{m}%
\sum_{i=1}^{n}M_{ik},  \label{2.7}
\end{equation}%
which is obviously coarser than (\ref{2.2}), but perhaps more useful for
applications.
\end{remark}

The case of inner product spaces, in which we may provide a simpler
condition of equality, is of interest in applications.

\begin{theorem}
\label{t2.2}Let $\left( X,\left\Vert \cdot \right\Vert \right) $ be an inner
product space over the real or complex number field $\mathbb{K}$, $e_{k},$ $%
x_{i}\in H\backslash \left\{ 0\right\} ,$ $k\in \left\{ 1,\dots ,m\right\} ,$
$i\in \left\{ 1,\dots ,n\right\} .$ If $M_{ik}\geq 0$ for $i\in \left\{
1,\dots ,n\right\} ,$ $\left\{ 1,\dots ,n\right\} $ such that%
\begin{equation}
\left\Vert x_{i}\right\Vert -\func{Re}\left\langle x_{i},e_{k}\right\rangle
\leq M_{ik}\text{ \ for each \ }i\in \left\{ 1,\dots ,n\right\} ,\ k\in
\left\{ 1,\dots ,m\right\} ,  \label{2.8}
\end{equation}%
then we have the inequality%
\begin{equation}
\sum_{i=1}^{n}\left\Vert x_{i}\right\Vert \leq \left\Vert \frac{1}{m}%
\sum_{k=1}^{m}e_{k}\right\Vert \left\Vert \sum_{i=1}^{n}x_{i}\right\Vert +%
\frac{1}{m}\sum_{k=1}^{m}\sum_{i=1}^{n}M_{ik}.  \label{2.9}
\end{equation}%
The case of equality holds in (\ref{2.9}) if and only if%
\begin{equation}
\sum_{i=1}^{n}\left\Vert x_{i}\right\Vert \geq \frac{1}{m}%
\sum_{k=1}^{m}\sum_{i=1}^{n}M_{ik}  \label{2.9.a}
\end{equation}%
and%
\begin{equation}
\sum_{i=1}^{n}x_{i}=\frac{m\left( \sum_{i=1}^{n}\left\Vert x_{i}\right\Vert -%
\frac{1}{m}\sum_{k=1}^{m}\sum_{i=1}^{n}M_{ik}\right) }{\left\Vert
\sum_{k=1}^{m}e_{k}\right\Vert ^{2}}\sum_{k=1}^{m}e_{k}.  \label{2.9.1}
\end{equation}
\end{theorem}

\begin{proof}
As in the proof of Theorem \ref{t2.1}, we have%
\begin{equation}
\sum_{i=1}^{n}\left\Vert x_{i}\right\Vert \leq \func{Re}\left\langle \frac{1%
}{m}\sum_{k=1}^{m}e_{k},\sum_{i=1}^{n}x_{i}\right\rangle +\frac{1}{m}%
\sum_{k=1}^{m}\sum_{i=1}^{n}M_{ik},  \label{2.11}
\end{equation}%
and $\sum_{k=1}^{m}e_{k}\neq 0.$

On utilising the Schwarz inequality in the inner product space $\left(
H;\left\langle \cdot ,\cdot \right\rangle \right) $ for $%
\sum_{i=1}^{n}x_{i}, $ $\sum_{k=1}^{m}e_{k},$ we have%
\begin{align}
\left\Vert \sum_{i=1}^{n}x_{i}\right\Vert \left\Vert
\sum_{k=1}^{m}e_{k}\right\Vert & \geq \left\vert \left\langle
\sum_{i=1}^{n}x_{i},\sum_{k=1}^{m}e_{k}\right\rangle \right\vert \geq
\left\vert \func{Re}\left\langle
\sum_{i=1}^{n}x_{i},\sum_{k=1}^{m}e_{k}\right\rangle \right\vert
\label{2.12} \\
& \geq \func{Re}\left\langle
\sum_{i=1}^{n}x_{i},\sum_{k=1}^{m}e_{k}\right\rangle .  \notag
\end{align}%
By (\ref{2.11}) and (\ref{2.12}) we deduce (\ref{2.9}).

Taking the norm in (\ref{2.9.1}) and using (\ref{2.9.a}), we have%
\begin{equation*}
\left\Vert \sum_{i=1}^{n}x_{i}\right\Vert =\frac{m\left(
\sum_{i=1}^{n}\left\Vert x_{i}\right\Vert -\frac{1}{m}\sum_{k=1}^{m}%
\sum_{i=1}^{n}M_{ik}\right) }{\left\Vert \sum_{k=1}^{m}e_{k}\right\Vert },
\end{equation*}%
showing that the equality holds in (\ref{2.9}).

Conversely, if the case of equality holds in (\ref{2.9}), then it must hold
in all the inequalities used to prove (\ref{2.9}). Therefore we have%
\begin{equation}
\left\Vert x_{i}\right\Vert =\func{Re}\left\langle x_{i},e_{k}\right\rangle
+M_{ik}\text{ \ for each \ }i\in \left\{ 1,\dots ,n\right\} ,\ k\in \left\{
1,\dots ,m\right\} ,  \label{2.13}
\end{equation}%
\begin{equation}
\left\Vert \sum_{i=1}^{n}x_{i}\right\Vert \left\Vert
\sum_{k=1}^{m}e_{k}\right\Vert =\left\vert \left\langle
\sum_{i=1}^{n}x_{i},\sum_{k=1}^{m}e_{k}\right\rangle \right\vert
\label{2.14}
\end{equation}%
and%
\begin{equation}
\func{Im}\left\langle \sum_{i=1}^{n}x_{i},\sum_{k=1}^{m}e_{k}\right\rangle
=0.  \label{2.15}
\end{equation}%
From (\ref{2.13}), on summing over $i$ and $k,$ we get%
\begin{equation}
\func{Re}\left\langle \sum_{i=1}^{n}x_{i},\sum_{k=1}^{m}e_{k}\right\rangle
=m\sum_{i=1}^{n}\left\Vert x_{i}\right\Vert
-\sum_{k=1}^{m}\sum_{i=1}^{n}M_{ik}.  \label{2.16}
\end{equation}%
On the other hand, by the use of the following identity in inner product
spaces,%
\begin{equation*}
\left\Vert u-\frac{\left\langle u,v\right\rangle v}{\left\Vert v\right\Vert
^{2}}\right\Vert ^{2}=\frac{\left\Vert u\right\Vert ^{2}\left\Vert
v\right\Vert ^{2}-\left\vert \left\langle u,v\right\rangle \right\vert ^{2}}{%
\left\Vert v\right\Vert ^{2}},\quad v\neq 0;
\end{equation*}%
the relation (\ref{2.14}) holds if and only if%
\begin{equation*}
\sum_{i=1}^{n}x_{i}=\frac{\left\langle
\sum_{i=1}^{n}x_{i},\sum_{k=1}^{m}e_{k}\right\rangle }{\left\Vert
\sum_{k=1}^{m}e_{k}\right\Vert ^{2}}\sum_{k=1}^{m}e_{k},
\end{equation*}%
giving, from (\ref{2.15}) and (\ref{2.16}), that%
\begin{equation*}
\sum_{i=1}^{n}x_{i}=\frac{m\sum_{i=1}^{n}\left\Vert x_{i}\right\Vert
-\sum_{k=1}^{m}\sum_{i=1}^{n}M_{ik}}{\left\Vert
\sum_{k=1}^{m}e_{k}\right\Vert ^{2}}\sum_{k=1}^{m}e_{k}.
\end{equation*}%
If the inequality holds in (\ref{2.9}), then obviously (\ref{2.9.a}) is
valid, and the theorem is proved.
\end{proof}

\begin{remark}
\label{r2.2}If in the above theorem the vectors $\left\{ e_{k}\right\} _{k=%
\overline{1,m}}$ are assumed to be orthogonal, then (\ref{2.9}) becomes:%
\begin{equation}
\sum_{i=1}^{n}\left\Vert x_{i}\right\Vert \leq \frac{1}{m}\left(
\sum_{k=1}^{m}\left\Vert e_{k}\right\Vert ^{2}\right) ^{\frac{1}{2}%
}\left\Vert \sum_{i=1}^{n}x_{i}\right\Vert +\frac{1}{m}\sum_{k=1}^{m}%
\sum_{i=1}^{n}M_{ik}.  \label{2.16.a}
\end{equation}%
Moreover, if $\left\{ e_{k}\right\} _{k=\overline{1,m}}$ is an orthonormal
family, then (\ref{2.16.a}) becomes 
\begin{equation}
\sum_{i=1}^{n}\left\Vert x_{i}\right\Vert \leq \frac{\sqrt{m}}{m}\left\Vert
\sum_{i=1}^{n}x_{i}\right\Vert +\frac{1}{m}\sum_{k=1}^{m}%
\sum_{i=1}^{n}M_{ik},  \label{2.16.b}
\end{equation}%
which has been obtained in \cite{SSD3}.
\end{remark}

Before we provide some natural consequences of Theorem \ref{t2.2}, we need
some preliminary results concerning reverses of Schwarz's inequality in
inner product spaces (see for instance \cite[p. 27]{SSD2}).

\begin{lemma}
\label{l2.1}Let $\left( X,\left\Vert \cdot \right\Vert \right) $ be an inner
product space over the real or complex number field $\mathbb{K}$ and $x,a\in
H,$ $r>0.$ If $\left\Vert x-a\right\Vert \leq r,$ then we have the inequality%
\begin{equation}
\left\Vert x\right\Vert \left\Vert a\right\Vert -\func{Re}\left\langle
x,a\right\rangle \leq \frac{1}{2}r^{2}.  \label{2.17}
\end{equation}%
The case of equality holds in (\ref{2.17}) if and only if%
\begin{equation}
\left\Vert x-a\right\Vert =r\text{ \ and \ }\left\Vert x\right\Vert
=\left\Vert a\right\Vert .  \label{2.18}
\end{equation}
\end{lemma}

\begin{proof}
The condition $\left\Vert x-a\right\Vert \leq r$ is clearly equivalent to%
\begin{equation}
\left\Vert x\right\Vert ^{2}+\left\Vert a\right\Vert ^{2}\leq 2\func{Re}%
\left\langle x,a\right\rangle +r^{2}.  \label{2.19}
\end{equation}%
Since%
\begin{equation}
2\left\Vert x\right\Vert \left\Vert a\right\Vert \leq \left\Vert
x\right\Vert ^{2}+\left\Vert a\right\Vert ^{2},  \label{2.20}
\end{equation}%
with equality if and only if $\left\Vert x\right\Vert =\left\Vert
a\right\Vert ,$ hence by (\ref{2.19}) and (\ref{2.20}) we deduce (\ref{2.17}%
).

The case of equality is obvious.
\end{proof}

Utilising the above lemma we may state the following corollary of Theorem %
\ref{t2.2}.

\begin{corollary}
\label{c2.1}Let $\left( H;\left\langle \cdot ,\cdot \right\rangle \right) ,$ 
$e_{k},$ $x_{i}$ be as in Theorem \ref{t2.2}. If $r_{ik}>0,$ $i\in \left\{
1,\dots ,n\right\} ,$ $k\in \left\{ 1,\dots ,m\right\} $ such that%
\begin{equation}
\left\Vert x_{i}-e_{k}\right\Vert \leq r_{ik}\text{ \ for each \ }i\in
\left\{ 1,\dots ,n\right\} \text{ and }k\in \left\{ 1,\dots ,m\right\} ,
\label{2.21}
\end{equation}%
then we have the inequality 
\begin{equation}
\sum_{i=1}^{n}\left\Vert x_{i}\right\Vert \leq \left\Vert \frac{1}{m}%
\sum_{k=1}^{m}e_{k}\right\Vert \left\Vert \sum_{i=1}^{n}x_{i}\right\Vert +%
\frac{1}{2m}\sum_{k=1}^{m}\sum_{i=1}^{n}r_{ik}^{2}.  \label{2.22}
\end{equation}%
The equality holds in (\ref{2.22}) if and only if%
\begin{equation*}
\sum_{i=1}^{n}\left\Vert x_{i}\right\Vert \geq \frac{1}{2m}%
\sum_{k=1}^{m}\sum_{i=1}^{n}r_{ik}^{2}
\end{equation*}%
and%
\begin{equation*}
\sum_{i=1}^{n}x_{i}=\frac{m\left( \sum_{i=1}^{n}\left\Vert x_{i}\right\Vert -%
\frac{1}{2m}\sum_{k=1}^{m}\sum_{i=1}^{n}r_{ik}^{2}\right) }{\left\Vert
\sum_{k=1}^{m}e_{k}\right\Vert ^{2}}\sum_{k=1}^{m}e_{k}.
\end{equation*}
\end{corollary}

The following lemma may provide another sufficient condition for (\ref{2.8})
to hold (see also \cite[p. 28]{SSD2}).

\begin{lemma}
\label{l2.2}Let $\left( H;\left\langle \cdot ,\cdot \right\rangle \right) $
be an inner product space over the real or complex number field $\mathbb{K}$
and $x,y\in H,$ $M\geq m>0.$ If either%
\begin{equation}
\func{Re}\left\langle My-x,x-my\right\rangle \geq 0  \label{2.23}
\end{equation}%
or, equivalently,%
\begin{equation}
\left\Vert x-\frac{m+M}{2}y\right\Vert \leq \frac{1}{2}\left( M-m\right)
\left\Vert y\right\Vert ,  \label{2.24}
\end{equation}%
holds, then%
\begin{equation}
\left\Vert x\right\Vert \left\Vert y\right\Vert -\func{Re}\left\langle
x,y\right\rangle \leq \frac{1}{4}\cdot \frac{\left( M-m\right) ^{2}}{m+M}%
\left\Vert y\right\Vert ^{2}.  \label{2.25}
\end{equation}%
The case of equality holds in (\ref{2.25}) if and only if the equality case
is realised in (\ref{2.23}) and%
\begin{equation*}
\left\Vert x\right\Vert =\frac{M+m}{2}\left\Vert y\right\Vert .
\end{equation*}
\end{lemma}

The proof is obvious by Lemma \ref{l2.1} for $a=\frac{M+m}{2}y$ and $r=\frac{%
1}{2}\left( M-m\right) \left\Vert y\right\Vert .$

Finally, the following corollary of Theorem \ref{t2.2} may be stated.

\begin{corollary}
\label{c2.2}Assume that $\left( H,\left\langle \cdot ,\cdot \right\rangle
\right) ,$ $e_{k},$ $x_{i}$ are as in Theorem \ref{t2.2}. If $M_{ik}\geq
m_{ik}>0$ satisfy the condition%
\begin{equation*}
\func{Re}\left\langle M_{k}e_{k}-x_{i},x_{i}-\mu _{k}e_{k}\right\rangle \geq
0\text{ \ for each \ }i\in \left\{ 1,\dots ,n\right\} \text{ and }k\in
\left\{ 1,\dots ,m\right\} ,
\end{equation*}%
then%
\begin{equation*}
\sum_{i=1}^{n}\left\Vert x_{i}\right\Vert \leq \left\Vert \frac{1}{m}%
\sum_{k=1}^{m}e_{k}\right\Vert \left\Vert \sum_{i=1}^{n}x_{i}\right\Vert +%
\frac{1}{4m}\sum_{k=1}^{m}\sum_{i=1}^{n}\frac{\left( M_{ik}-m_{ik}\right)
^{2}}{M_{ik}+m_{ik}}\left\Vert e_{k}\right\Vert ^{2}.
\end{equation*}
\end{corollary}

\section{Applications for Complex Numbers}

Let $\mathbb{C}$ be the field of complex numbers. If $z=\func{Re}z+i\func{Im}%
z,$ then by $\left\vert \cdot \right\vert _{p}:\mathbb{C}\rightarrow \lbrack
0,\infty ),$ $p\in \left[ 1,\infty \right] $ we define the $p-$\textit{%
modulus }of $z$ as%
\begin{equation*}
\left\vert z\right\vert _{p}:=\left\{ 
\begin{array}{ll}
\max \left\{ \left\vert \func{Re}z\right\vert ,\left\vert \func{Im}%
z\right\vert \right\} & \text{if \ }p=\infty , \\ 
&  \\ 
\left( \left\vert \func{Re}z\right\vert ^{p}+\left\vert \func{Im}%
z\right\vert ^{p}\right) ^{\frac{1}{p}} & \text{if \ }p\in \lbrack 1,\infty
),%
\end{array}%
\right.
\end{equation*}%
where $\left\vert a\right\vert ,$ $a\in \mathbb{R}$ is the usual modulus of
the real number $a.$ Obviously, for $p=2,$ we recapture the usual modulus of
a complex number.

It is well known that $\left( \mathbb{C},\left\vert \cdot \right\vert
_{p}\right) ,$ $p\in \left[ 1,\infty \right] $ is a Banach space over the
complex number field $\mathbb{C}$.

Consider the Banach space $\left( \mathbb{C},\left\vert \cdot \right\vert
_{1}\right) $ and $F:\mathbb{C\rightarrow C}$, $F\left( z\right) =az$ with $%
a\in \mathbb{C}$, $a\neq 0.$ Obviously, $F$ is linear on $\mathbb{C}$. For $%
z\neq 0,$ we have%
\begin{equation*}
\frac{\left\vert F\left( z\right) \right\vert }{\left\vert z\right\vert _{1}}%
=\frac{\left\vert a\right\vert \left\vert z\right\vert }{\left\vert
z\right\vert _{1}}=\frac{\left\vert a\right\vert \sqrt{\left\vert \func{Re}%
z\right\vert ^{2}+\left\vert \func{Im}z\right\vert ^{2}}}{\left\vert \func{Re%
}z\right\vert +\left\vert \func{Im}z\right\vert }\leq \left\vert
a\right\vert .
\end{equation*}%
Since, for $z_{0}=1,$ we have $\left\vert F\left( z_{0}\right) \right\vert
=\left\vert a\right\vert $ and $\left\vert z_{0}\right\vert _{1}=1,$ hence%
\begin{equation*}
\left\Vert F\right\Vert _{1}:=\sup_{z\neq 0}\frac{\left\vert F\left(
z\right) \right\vert }{\left\vert z\right\vert _{1}}=\left\vert a\right\vert
,
\end{equation*}%
showing that $F$ is a bounded linear functional on $\left( \mathbb{C}%
,\left\vert \cdot \right\vert _{1}\right) $ and $\left\Vert F\right\Vert
_{1}=\left\vert a\right\vert .$

We can apply Theorem \ref{t2.1} to state the following reverse of the
generalised triangle inequality\ for complex numbers.

\begin{proposition}
\label{p4.1}Let $a_{k},$ $x_{j}\in \mathbb{C}$, $k\in \left\{ 1,\dots
,m\right\} $ and $j\in \left\{ 1,\dots ,n\right\} .$ If there exist the
constants $M_{jk}\geq 0,$ $k\in \left\{ 1,\dots ,m\right\} ,$ $j\in \left\{
1,\dots ,n\right\} $ such that%
\begin{equation}
\left\vert \func{Re}x_{j}\right\vert +\left\vert \func{Im}x_{j}\right\vert
\leq \func{Re}a_{k}\cdot \func{Re}x_{j}-\func{Im}a_{k}\cdot \func{Im}%
x_{j}+M_{jk}  \label{4.1}
\end{equation}%
for each $j\in \left\{ 1,\dots ,n\right\} $ and $k\in \left\{ 1,\dots
,m\right\} ,$ then%
\begin{multline}
\quad \sum_{j=1}^{n}\left[ \left\vert \func{Re}x_{j}\right\vert +\left\vert 
\func{Im}x_{j}\right\vert \right]  \label{4.2} \\
\leq \frac{1}{m}\left\vert \sum_{k=1}^{m}a_{k}\right\vert \left[ \left\vert
\sum_{j=1}^{n}\func{Re}x_{j}\right\vert +\left\vert \sum_{j=1}^{n}\func{Im}%
x_{j}\right\vert \right] +\frac{1}{m}\sum_{k=1}^{m}\sum_{j=1}^{n}M_{jk}.\quad
\end{multline}
\end{proposition}

The proof follows by Theorem \ref{t2.1} applied for the Banach space $\left( 
\mathbb{C},\left\vert \cdot \right\vert _{1}\right) $ and $F_{k}\left(
z\right) =a_{k}z,$ $k\in \left\{ 1,\dots ,m\right\} $ on taking into account
that:%
\begin{equation*}
\left\Vert \sum_{k=1}^{m}F_{k}\right\Vert _{1}=\left\vert
\sum_{k=1}^{m}a_{k}\right\vert .
\end{equation*}%
Now, consider the Banach space $\left( \mathbb{C},\left\vert \cdot
\right\vert _{\infty }\right) .$ If $F\left( z\right) =dz,$ then for $z\neq 0
$ we have%
\begin{equation*}
\frac{\left\vert F\left( z\right) \right\vert }{\left\vert z\right\vert
_{\infty }}=\frac{\left\vert d\right\vert \left\vert z\right\vert }{%
\left\vert z\right\vert _{\infty }}=\frac{\left\vert d\right\vert \sqrt{%
\left\vert \func{Re}z\right\vert ^{2}+\left\vert \func{Im}z\right\vert ^{2}}%
}{\max \left\{ \left\vert \func{Re}z\right\vert ,\left\vert \func{Im}%
z\right\vert \right\} }\leq \sqrt{2}\left\vert d\right\vert .
\end{equation*}%
Since, for $z_{0}=1+i,$ we have $\left\vert F\left( z_{0}\right) \right\vert
=\sqrt{2}\left\vert d\right\vert ,$ $\left\vert z_{0}\right\vert _{\infty
}=1,$ hence%
\begin{equation*}
\left\Vert F\right\Vert _{\infty }:=\sup_{z\neq 0}\frac{\left\vert F\left(
z\right) \right\vert }{\left\vert z\right\vert _{\infty }}=\sqrt{2}%
\left\vert d\right\vert ,
\end{equation*}%
showing that $F$ is a bounded linear functional on $\left( \mathbb{C}%
,\left\vert \cdot \right\vert _{\infty }\right) $ and $\left\Vert
F\right\Vert _{\infty }=\sqrt{2}\left\vert d\right\vert .$

If we apply Theorem \ref{t2.1}, then we can state the following reverse of
the generalised triangle inequality for complex numbers.

\begin{proposition}
\label{p4.2}Let $a_{k},$ $x_{j}\in \mathbb{C}$, $k\in \left\{ 1,\dots
,m\right\} $ and $j\in \left\{ 1,\dots ,n\right\} .$ If there exist the
constants $M_{jk}\geq 0,$ $k\in \left\{ 1,\dots ,m\right\} ,$ $j\in \left\{
1,\dots ,n\right\} $ such that%
\begin{equation*}
\max \left\{ \left\vert \func{Re}x_{j}\right\vert ,\left\vert \func{Im}%
x_{j}\right\vert \right\} \leq \func{Re}a_{k}\cdot \func{Re}x_{j}-\func{Im}%
a_{k}\cdot \func{Im}x_{j}+M_{jk}
\end{equation*}%
for each $j\in \left\{ 1,\dots ,n\right\} $ and $k\in \left\{ 1,\dots
,m\right\} ,$ then%
\begin{multline}
\sum_{j=1}^{n}\max \left\{ \left\vert \func{Re}x_{j}\right\vert ,\left\vert 
\func{Im}x_{j}\right\vert \right\}   \label{4.3} \\
\leq \frac{\sqrt{2}}{m}\left\vert \sum_{k=1}^{m}a_{k}\right\vert \max
\left\{ \left\vert \sum_{j=1}^{n}\func{Re}x_{j}\right\vert ,\left\vert
\sum_{j=1}^{n}\func{Im}x_{j}\right\vert \right\} +\frac{1}{m}%
\sum_{k=1}^{m}\sum_{j=1}^{n}M_{jk}.
\end{multline}
\end{proposition}

Finally, consider the Banach space $\left( \mathbb{C},\left\vert \cdot
\right\vert _{2p}\right) $ with $p\geq 1.$

Let $F:\mathbb{C\rightarrow C}$, $F\left( z\right) =cz.$ By H\"{o}lder's
inequality, we have%
\begin{equation*}
\frac{\left\vert F\left( z\right) \right\vert }{\left\vert z\right\vert _{2p}%
}=\frac{\left\vert c\right\vert \sqrt{\left\vert \func{Re}z\right\vert
^{2}+\left\vert \func{Im}z\right\vert ^{2}}}{\left( \left\vert \func{Re}%
z\right\vert ^{2p}+\left\vert \func{Im}z\right\vert ^{2p}\right) ^{\frac{1}{%
2p}}}\leq 2^{\frac{1}{2}-\frac{1}{2p}}\left\vert c\right\vert .
\end{equation*}%
Since, for $z_{0}=1+i$ we have $\left\vert F\left( z_{0}\right) \right\vert
=2^{\frac{1}{2}}\left\vert c\right\vert ,$ $\left\vert z_{0}\right\vert =2^{%
\frac{1}{2p}}$ $\left( p\geq 1\right) ,$ hence%
\begin{equation*}
\left\Vert F\right\Vert _{2p}:=\sup_{z\neq 0}\frac{\left\vert F\left(
z\right) \right\vert }{\left\vert z\right\vert _{2p}}=2^{\frac{1}{2}-\frac{1%
}{2p}}\left\vert c\right\vert ,
\end{equation*}%
showing that $F$ is a bounded linear functional on $\left( \mathbb{C}%
,\left\vert \cdot \right\vert _{2p}\right) ,$ $p\geq 1$ and $\left\Vert
F\right\Vert _{2p}=2^{\frac{1}{2}-\frac{1}{2p}}\left\vert c\right\vert .$

If we apply Theorem \ref{t2.1}, then we can state the following proposition.

\begin{proposition}
\label{p4.3}Let $a_{k},$ $x_{j},$ $M_{jk}$ be as in Proposition \ref{p4.2}.
If%
\begin{equation*}
\left[ \left\vert \func{Re}x_{j}\right\vert ^{2p}+\left\vert \func{Im}%
x_{j}\right\vert ^{2p}\right] ^{\frac{1}{2p}}\leq \func{Re}a_{k}\cdot \func{%
Re}x_{j}-\func{Im}a_{k}\cdot \func{Im}x_{j}+M_{jk}
\end{equation*}%
for each $j\in \left\{ 1,\dots ,n\right\} $ and $k\in \left\{ 1,\dots
,m\right\} ,$ then%
\begin{multline}
\sum_{j=1}^{n}\left[ \left\vert \func{Re}x_{j}\right\vert ^{2p}+\left\vert 
\func{Im}x_{j}\right\vert ^{2p}\right] ^{\frac{1}{2p}}  \label{4.4} \\
\leq \frac{2^{\frac{1}{2}-\frac{1}{2p}}}{m}\left\vert
\sum_{k=1}^{m}a_{k}\right\vert \left[ \left\vert \sum_{j=1}^{n}\func{Re}%
x_{j}\right\vert ^{2p}+\left\vert \sum_{j=1}^{n}\func{Im}x_{j}\right\vert
^{2p}\right] ^{\frac{1}{2p}}+\frac{1}{m}\sum_{k=1}^{m}\sum_{j=1}^{n}M_{jk}.
\end{multline}%
where $p\geq 1.$
\end{proposition}

\begin{remark}
If in the above proposition we choose $p=1,$ then we have the following
reverse of the generalised triangle inequality for complex numbers%
\begin{equation*}
\sum_{j=1}^{n}\left\vert x_{j}\right\vert \leq \left\vert \frac{1}{m}%
\sum_{k=1}^{m}a_{k}\right\vert \left\vert \sum_{j=1}^{n}x_{j}\right\vert +%
\frac{1}{m}\sum_{k=1}^{m}\sum_{j=1}^{n}M_{jk}
\end{equation*}%
provided $x_{j},a_{k},$ $j\in \left\{ 1,\dots ,n\right\} $, $k\in \left\{
1,\dots ,m\right\} $ satisfy the assumption%
\begin{equation*}
\left\vert x_{j}\right\vert \leq \func{Re}a_{k}\cdot \func{Re}x_{j}-\func{Im}%
a_{k}\cdot \func{Im}x_{j}+M_{jk}
\end{equation*}%
for each $j\in \left\{ 1,\dots ,n\right\} $, $k\in \left\{ 1,\dots
,m\right\} .$ Here $\left\vert \cdot \right\vert $ is the usual modulus of a
complex number and $M_{jk}>0,j\in \left\{ 1,\dots ,n\right\} $, $k\in
\left\{ 1,\dots ,m\right\} $ are given.
\end{remark}

\end{document}